\documentclass{ws-procs9x6}            

\usepackage{cite}
\renewcommand{\cite}[1]{[\citen{#1}]}

\newcommand{\uvec}[2][3]{\boldsymbol{#2\mkern-#1mu}\mkern#1mu}

\newcommand*\xbar[1]{%
	\hbox{%
		\vbox{%
			\hrule height 0.5pt 
			\kern0.4ex
			\hbox{%
				\kern-0.05em
				\ensuremath{#1}%
				\kern-0.00em
			}%
		}%
	}%
}

\allowdisplaybreaks[1]

\newcommand{\jph}{{j+\frac{1}{2}}}
\newcommand{\jmh}{{j-\frac{1}{2}}}

\newcommand{\dx}{\Delta x}

\newcommand{\dt}{\Delta t}

\newcommand{\ves}{\varepsilon}

\newcommand{\Vbar}{\overline{\uvec{V}}}

\begin{document}
\title{New Smoothness Indicator Within an Active Flux Framework}

\author{Alina Chertock}
\address{Department of Mathematics, North Carolina State University,\\Raleigh, NC 27695, USA,\\acherto@ncsu.edu}

\author{Alexander Kurganov}
\address{Department of Mathematics and SUSTech International Center for Mathematics, Southern University of Science and Technology,
Shenzhen, 518055, China,\\alexander@sustech.edu.cn}

\author{Lorenzo Micalizzi}
\address{Department of Mathematics, North Carolina State University,\\Raleigh, NC 27695, USA,\\lmicali@ncsu.edu}

\begin{abstract}
In this work, we introduce a new smoothness indicator (SI), which is capable of detecting ``rough'' parts of the solutions computed by
active flux (AF) methods for hyperbolic (systems of) conservation laws. The new SI is based on measuring the difference between the two sets
of solutions (either cell averages and point values or cell averages on overlapping grids) evolved at each time step of AF methods. The key
idea in the derivation of the new SI is that in the ``rough'' parts of the evolved solutions, the difference is ${\cal O}(1)$, while
in the smooth areas, it is proportional to the order of the underlying AF method.

The performance of the new SI, that is, its ability to automatically and robustly detect ``rough'' parts of the computed solutions, is
illustrated on several numerical examples, in which the one-dimensional Euler equations of gas dynamics are numerically solved by a
recently introduced semi-discrete finite-volume AF method on overlapping grids.
\end{abstract}


\bodymatter
\section{Introduction}
We consider the one-dimensional (1-D) hyperbolic (system of) conservation laws 
\begin{equation}
\bm U_x+\bm F(\bm U)_x=\bm0,
\label{1}
\end{equation}
where $x$ is a spatial variable, $t$ is time, $\bm U\in\mathbb R^M$ is the vector of conservative variables, $\bm F\in\mathbb R^M$ is the
flux with real-diagonalizable Jacobian $\frac{\partial\uvec{F}}{\partial\uvec{U}}(\uvec{U})$.

In the past decade, active flux (AF) numerical methods for \eqref{1} have been developed and successfully applied to several hyperbolic
systems; we refer the reader to, e.g.,
\cite{eymann2013multidimensional,Zeng14,barsukow2021active,chudzik2021cartesian,Roe21,abgrall2023combination,abgrall2023extensions,AL24} and
references therein. The key idea of AF methods is to evolve several pieces of information on the computed solution, say, both cell averages
and point values at cell interfaces, or cell averages on overlapping grids as was done in our recent semi-discrete finite-volume (FV) AF
method introduced in \cite{ACKM}. This gives additional degrees of freedom, which can be used to enhance the accuracy of the resulting
scheme and successfully hybridize conservative and nonconservative numerical methods as some pieces of information (say, the point values)
can be accurately evolved using nonconservative (primitive variable) formulation of \eqref{1} without risking convergence to nonphyscial
weak solutions.

The main goal of this paper is to take another advantage of the two available sets of information and to design a new smoothness indicator
(SI), which is capable of detecting ``rough'' parts of the solutions computed by AF methods. The key idea in the derivation of the new SI is
that in the ``rough'' parts of the evolved solutions, the difference between the two solution versions is ${\cal O}(1)$, while in the smooth
areas, it is proportional to the order of the underlying AF method. We test the behavior of the SI on several numerical examples, in which
the one-dimensional (1-D) Euler equations of gas dynamics are numerically solved by the semi-discrete FV AF method from \cite{ACKM}. The
proposed SI can be used to develop several adaption strategies, which we will explore in our future work.

\section{Semi-Discrete FV AF Method---A Brief Overview}
We consider overlapping FV meshes consisting of uniform cells $I_j=[x_\jmh,x_\jph]$, $j=1,\ldots,N$ and $I_\jph=[x_j,x_{j+1}]$,
$j=0,\ldots,N$ with $x_{j+1}=x_\jph+\dx/2=x_j+\dx$. We assume that the computed cell averages
$\,\xbar{\bm U}_j(t):\approx\int_{I_j}\bm U(x,t)\,{\rm d}x$ are available at a certain time level $t$.

In addition, we consider an equivalent (for smooth solutions) nonconservative formulation of \eqref{1}:
\begin{equation}
\bm V_t+\widetilde{\bm F}(\bm V)_x=B(\bm V)\bm V_x,
\label{2}
\end{equation}
where $\bm V\in\mathbb R^M$ is the vector of primitive variables, $\widetilde{\bm F}:\mathbb R^M\to\mathbb R^M$, and
$B\in\mathbb R^{M\times M}$. We assume that the computed cell averages of $\bm V$ over the staggered cells,
$\,\xbar{\bm V}_\jph(t):\approx\int_{I_\jph}\bm V(x,t)\,{\rm d}x$, are also available at time $t$.

In the second-order semi-discrete FV AF method from \cite{ACKM}, the solution is evolved in time by numerically solving the following
system of ODEs:
\begin{equation}
\begin{aligned}
\frac{\rm d}{{\rm d}t}\,\xbar{\bm U}_j&=-\frac{1}{\dx}\Big[\bm{{\cal F}}_\jph-\bm{{\cal F}}_\jmh\Big],\\
\frac{\rm d}{{\rm d}t}\Vbar_\jph&=-\frac{1}{\dx}\bigg[\widetilde{\bm{{\cal F}}}_{j+1}-\widetilde{\bm{{\cal F}}}_j-\bm B_\jph\\
&\hspace*{1.2cm}-\frac{a^+_j}{a^+_j-a^-_j}\bm B_{\bm\Psi,j}+\frac{a^-_{j+1}}{a^+_{j+1}-a^-_{j+1}}\bm B_{\bm\Psi,j+1}\bigg],
\end{aligned}
\label{3}
\end{equation}
where $\bm{{\cal F}}_\jph$ are very simple numerical fluxes given by $\bm{{\cal F}}_\jph=\bm F\big(\bm U_\jph\big)$ with
$\bm U_\jph:=\bm U\big(\xbar{\bm V}_\jph\big)$ representing the point values of the conserved variables at cell interfaces. The
nonconservative system \eqref{2} was discretized using the path-conservative central-upwind (PCCU) scheme from \cite{diaz2019path}, and the
PCCU structures $\widetilde{\bm{{\cal F}}}_j$, $\bm B_\jph$, $\bm B_{\bm\Psi,j}$, $a^+_j$, and $a^-_j$ are defined in \cite{ACKM}. In the
numerical experiments reported in \S\ref{sec4}, we have integrated the ODE system \eqref{3} using the three-stage third-order strong
stability-preserving Runge-Kutta method (see, e.g., \cite{gottlieb2001strong}) with CFL number $0.25$. We note that the indexed quantities
in \eqref{3} and below are time-dependent, but we omit this dependence for the sake of brevity.

We would like to stress that after completing each time evolution step, the conservative cell averages at the new time level,
$\,\xbar{\bm U}_j(t+\dt)$, will likely be oscillatory as no limiting procedure was employed in the computation of the numerical fluxes 
$\bm{{\cal F}}_\jph$. At the same time, the oscillation-free cell averages $\Vbar_\jph$ are unreliable as solving the nonconservative system
\eqref{2} is known to lead nonphysical weak solutions; see, e.g., \cite{abgrall2010comment}. We therefore apply a conservative
post-processing that couples these two sets of cell averages and hence removes the oscillations in $\bm U$ and enforces $\bm V$ to converge
to the physically relevant solution; see \cite{ACKM} for details.

\section{New Smoothness Indicator (SI)}
In this section, we describe how to design a new SI, which can be computed after completing the time evolution step, but before the
post-processing. To this end, we introduce the following quantities:
\begin{equation}
\ves_j:=|\alpha(\,\xbar{\bm U}_j)-\alpha(\bm U(\bm V_j))|,\quad\bm V_j:=\frac{1}{2}\big(\,\xbar{\bm V}_\jmh+\xbar{\bm V}_\jph\big),
\label{4}
\end{equation}
where $\alpha$ is either one of the components of $\bm U$ or a function of $\bm U$. We then implement the noise-filtering from
\cite{chu2025new} and define
\begin{equation}
\widehat\ves_j:=\frac{1}{6}(\ves_{j-1}+4\ves_j+\ves_{j+1}).
\label{5}
\end{equation}
The SI \eqref{4}--\eqref{5} can be used to identify ``rough'' parts of the computed solution, for instance, as follows: We say that the
solution is smooth in cell $I_j$ if $\widehat\ves_j<K\widehat\ves_{\rm ave}$, where $\widehat\ves_{\rm ave}$ is the average of all
$\widehat\ves_j$ and $K$ is a positive tunable constant. This strategy is based on the following prediction: We expect the SI values to be
proportional to $(\dx)^2$ in smooth regions, coherently with the order of accuracy of the scheme; on contrary, in ``rough'' areas, we expect
$\widehat\ves_j$ to be ${\cal O}(1)$.

\section{Numerical Results}\label{sec4}
In this section, we investigate the behavior of the proposed SI on the 1-D Euler equations of gas dynamics. Their conservative (in terms of
the density, momentum and total energy) and primitive (in terms of the density, velocity and pressure) formulations can be found in
\cite{ACKM}. We use the same benchmarks as in Examples 1, 2, and 3 in \cite{chu2025new}, to which the reader is referred for the information
concerning the problem setups. In all of the numerical examples below, the quantity $\alpha$ in \eqref{4} is taken to be the momentum.

All of the solutions are computed by the second-order semi-discrete FV AF method from \cite{ACKM} and plotted along with the corresponding
SI values at the last time-step. In all SI plots below, we also plot two horizontal lines: one corresponding to $K\widehat\ves_{\rm ave}$
represents a borderline between the detected ``rough'' and smooth areas (for a given value of $K$) and the other one corresponding to
$C(\dx)^2$ is shown as a reference for the decay of the SI in smooth areas.

\begin{figure}[ht!]
\centerline{\includegraphics[trim=0.2cm 0.3cm 0.2cm 0.3cm, clip, width=5.6cm]{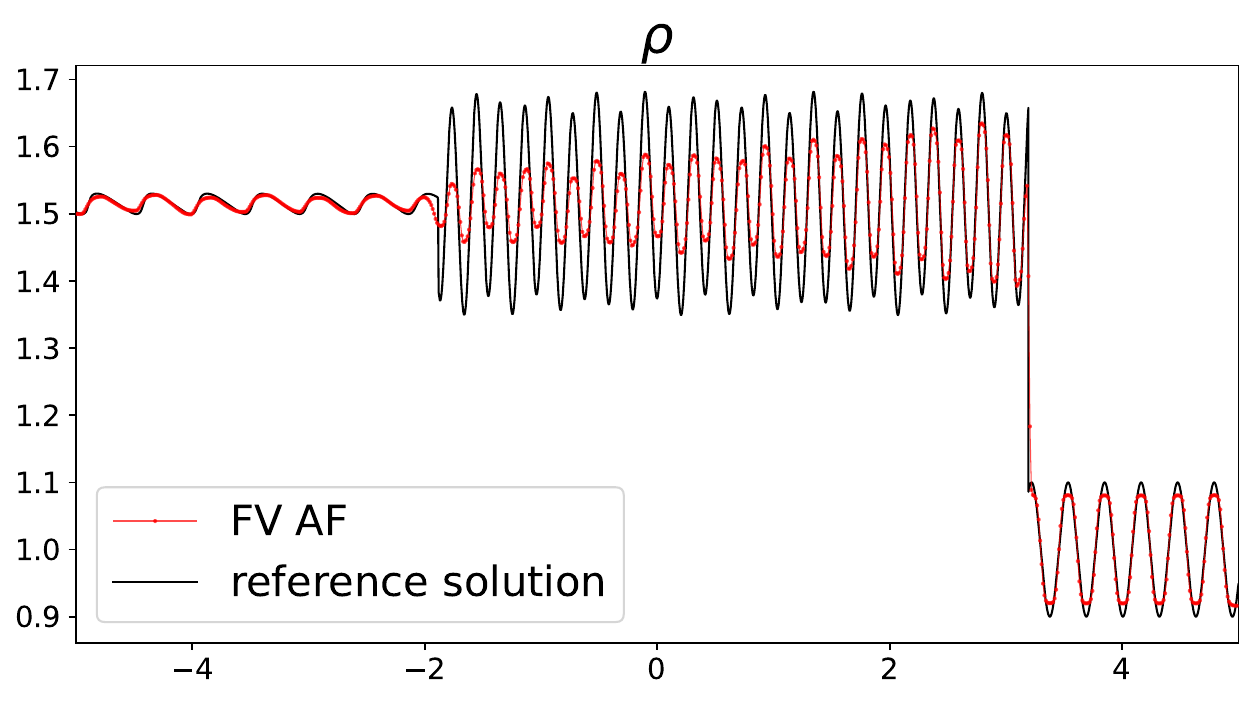}}
\vskip5pt
\centerline{\includegraphics[trim=0.3cm 0.3cm 0.2cm 0.2cm, clip, width=5.4cm]{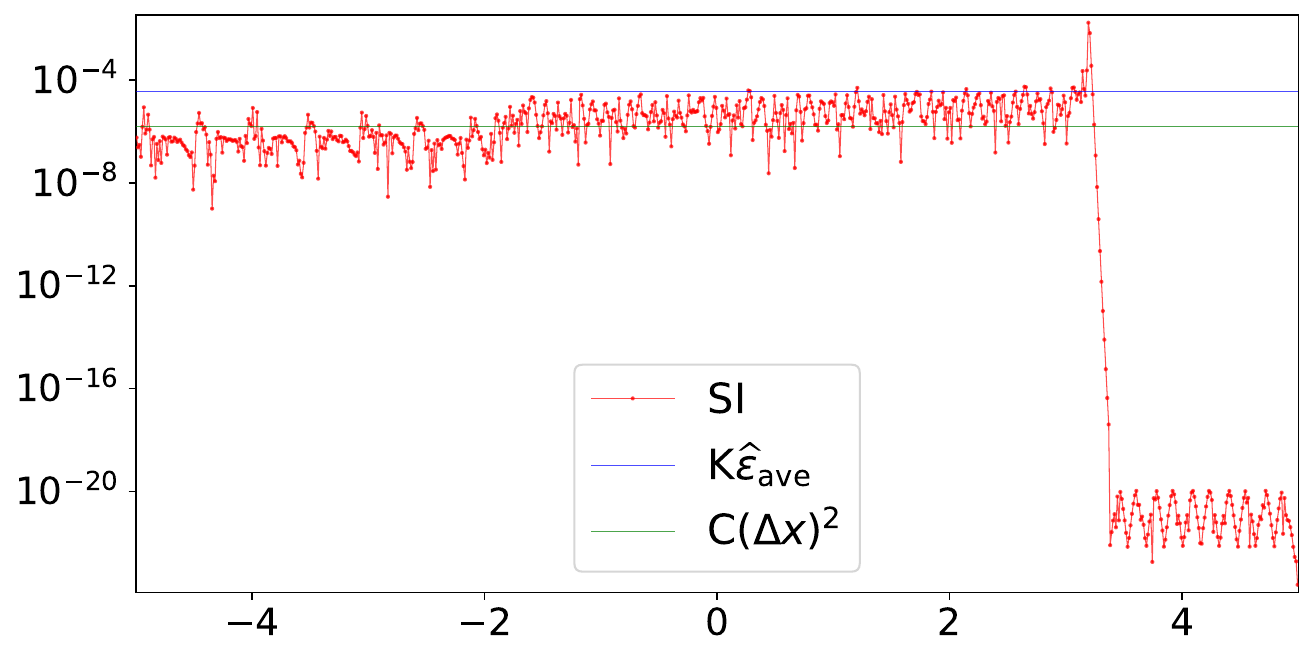}\hspace*{0.5cm}
	    \includegraphics[trim=0.3cm 0.3cm 0.2cm 0.2cm, clip, width=5.4cm]{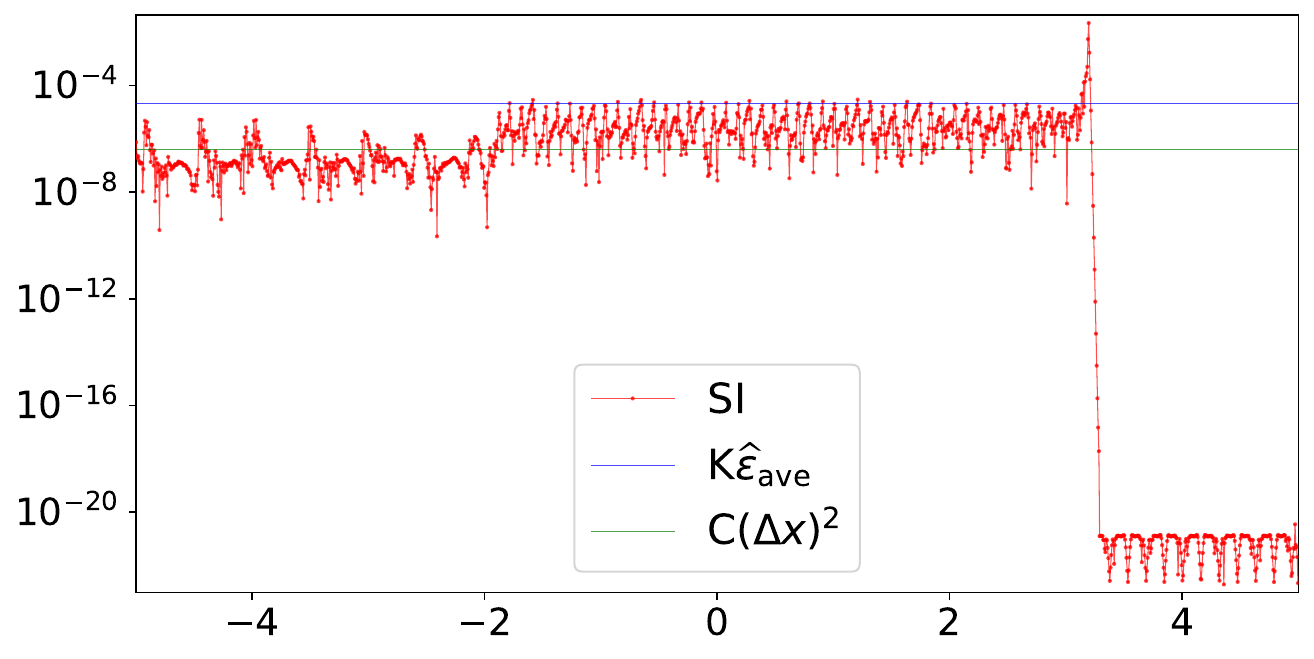}}
\vskip5pt
\centerline{\includegraphics[trim=0.3cm 0.3cm 0.2cm 0.2cm, clip, width=5.4cm]{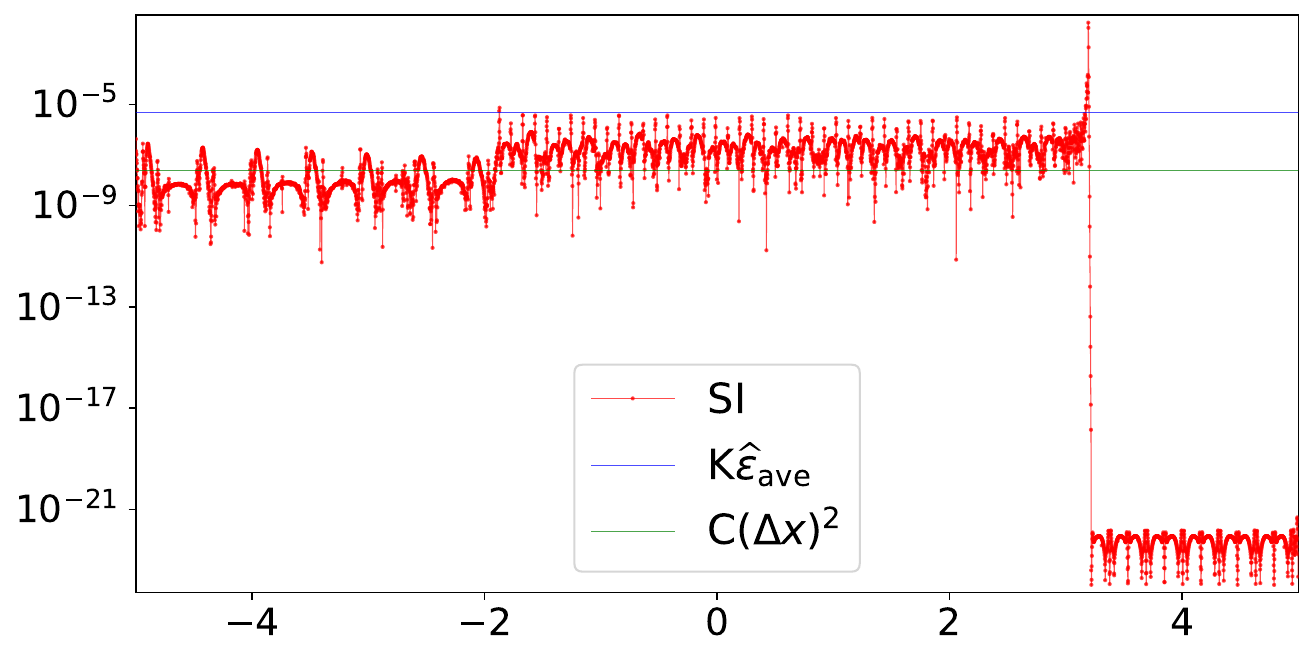}\hspace*{0.5cm}
            \includegraphics[trim=0.3cm 0.3cm 0.2cm 0.2cm, clip, width=5.4cm]{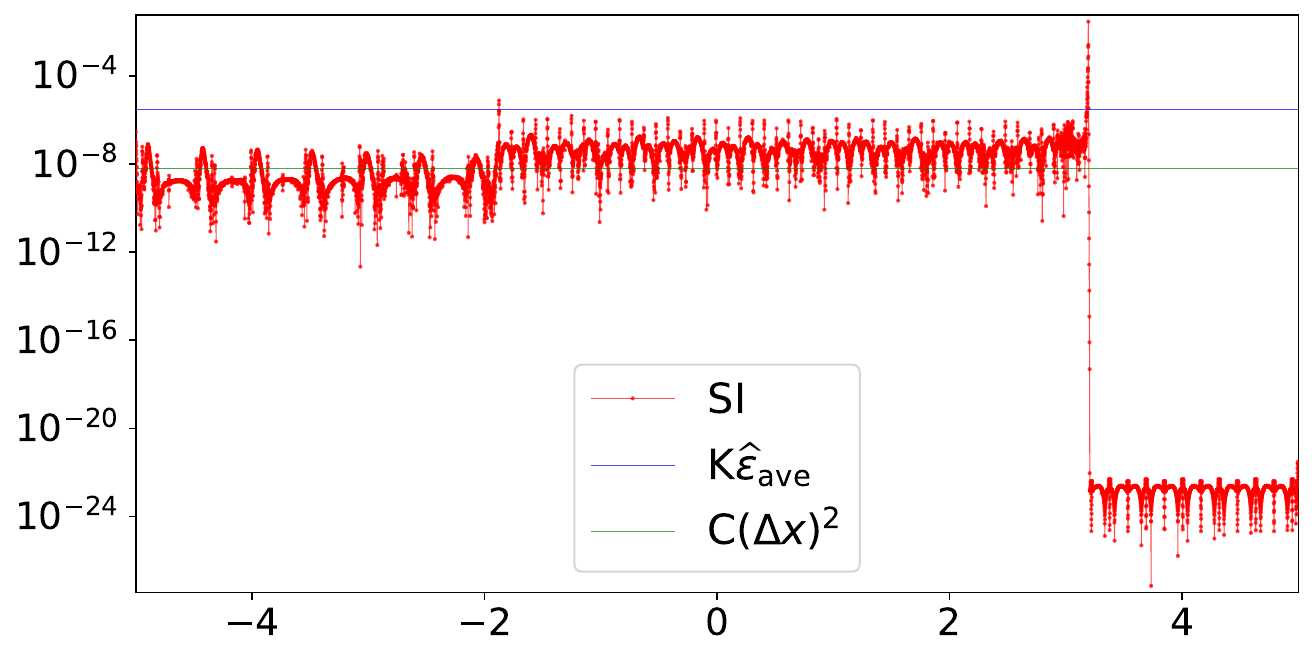}}
\caption{Example 1: Density computed using $N=800$ (top panel) and the corresponding SI values obtained on several meshes with $N=800$
(middle left), $N=1600$ (middle right), $N=6400$ (low left), and $N=12800$ (low right). $C=0.01$, $K=1$.\label{fig1}}
\end{figure}

\vskip5pt
\noindent
{\bf Example 1---shock-entropy wave interaction problem.}
The density computed with $N=800$ is plotted in Fig. \ref{fig1} (top panel). The corresponding SI values obtained on several mesh
refinements are shown in Fig. \ref{fig1} (middle and low rows). As one can clearly see, the shocks are detected (to detect a smaller one the
mesh had to be refined) and in the smooth region the SI decay is roughly $\sim(\dx)^2$ as expected.

\begin{figure}[ht!]
\centerline{\includegraphics[trim=0.2cm 0.3cm 0.2cm 0.3cm, clip, width=5.6cm]{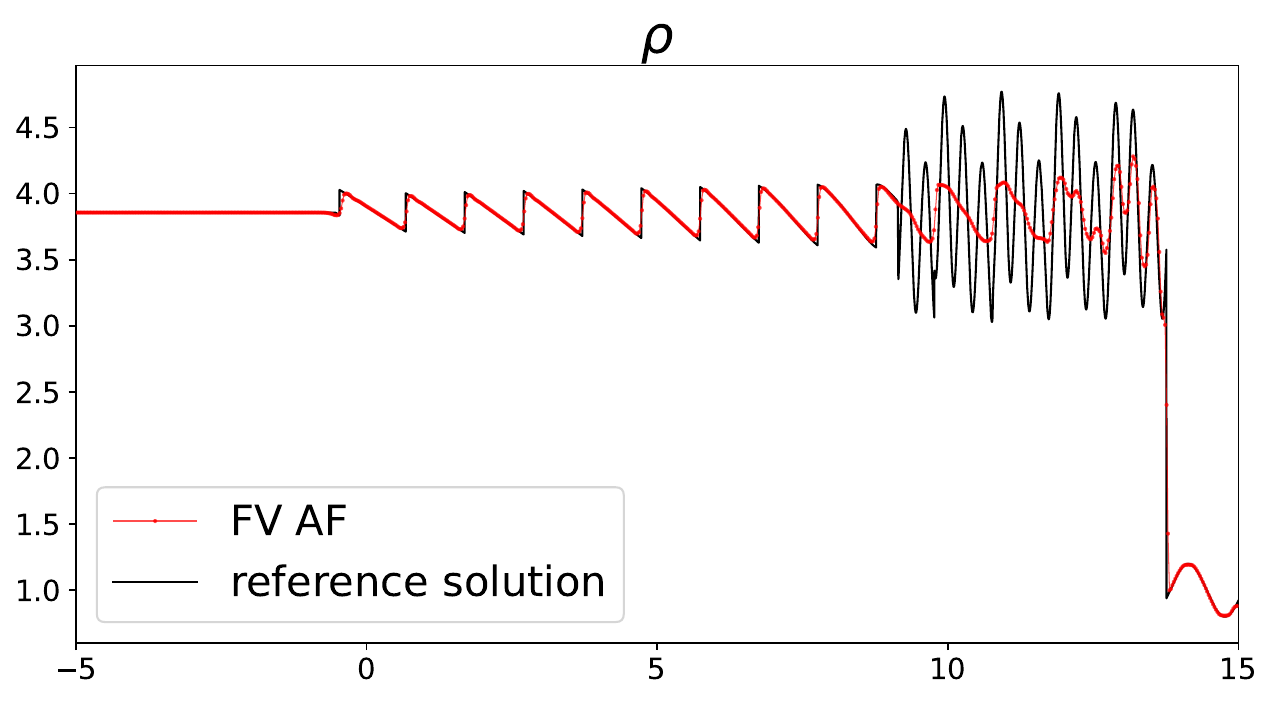}}
\vskip5pt
\centerline{\includegraphics[trim=0.3cm 0.3cm 0.2cm 0.2cm, clip, width=5.4cm]{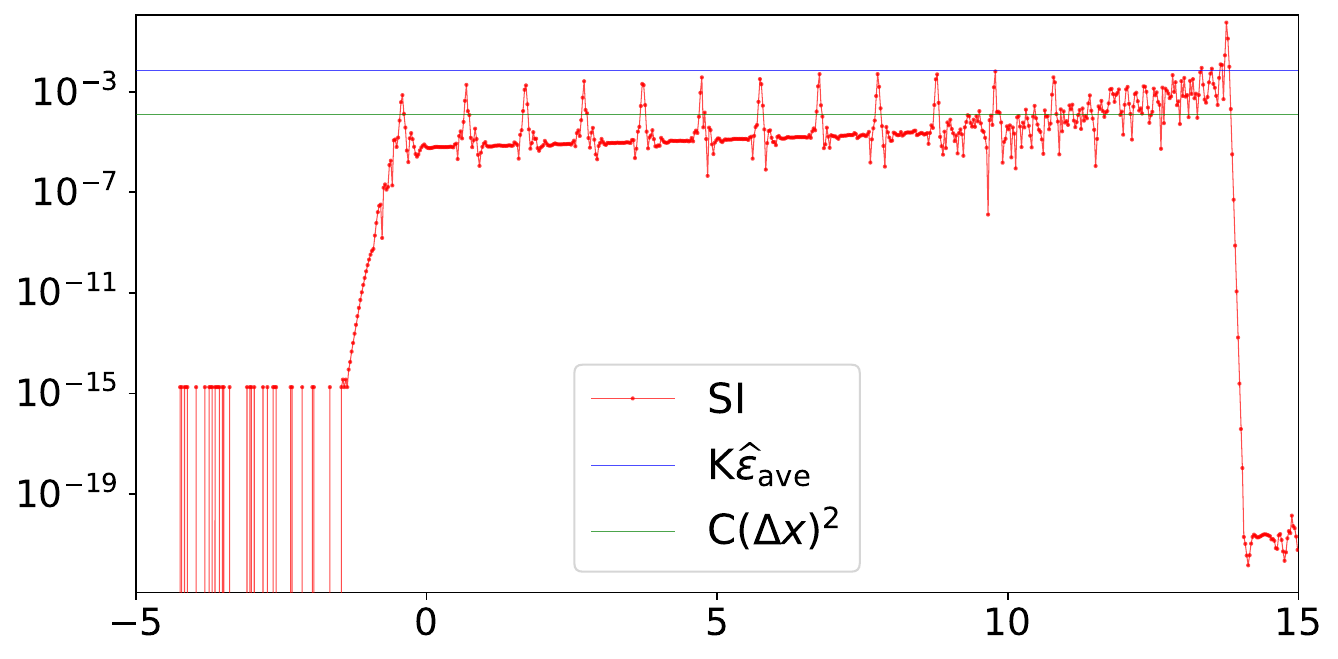}\hspace*{0.5cm}
            \includegraphics[trim=0.3cm 0.3cm 0.2cm 0.2cm, clip, width=5.4cm]{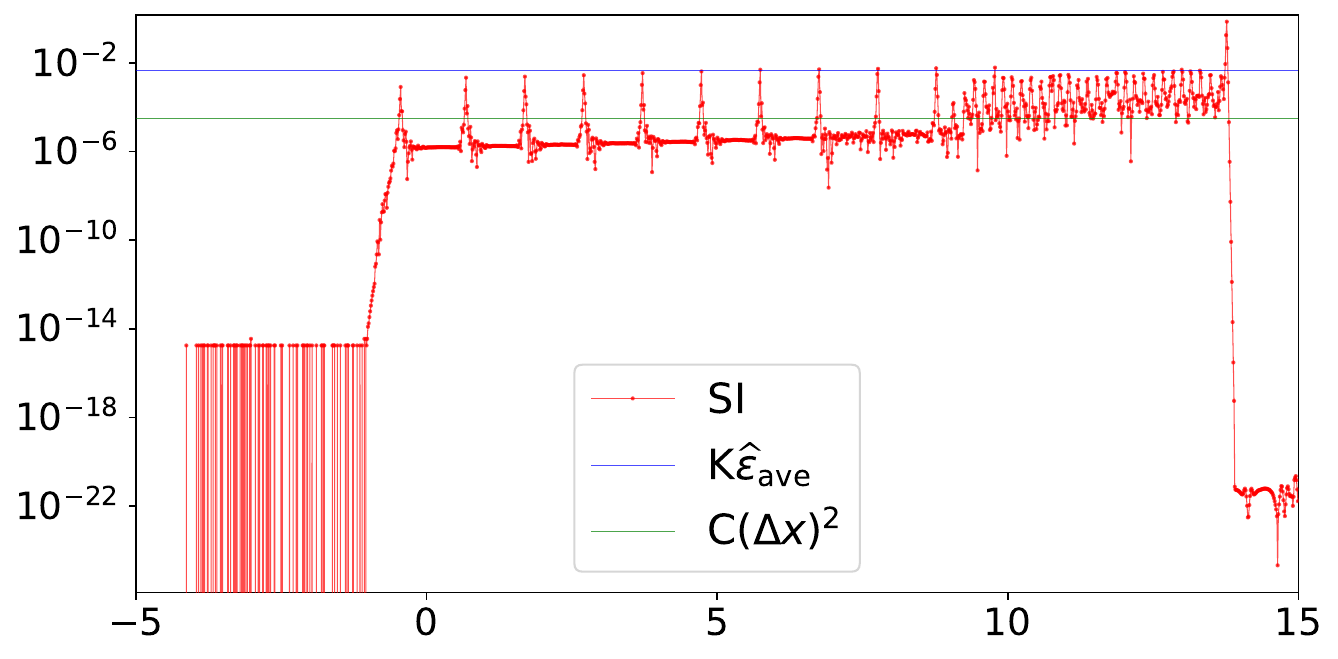}}
\vskip5pt
\centerline{\includegraphics[trim=0.3cm 0.3cm 0.2cm 0.2cm, clip, width=5.4cm]{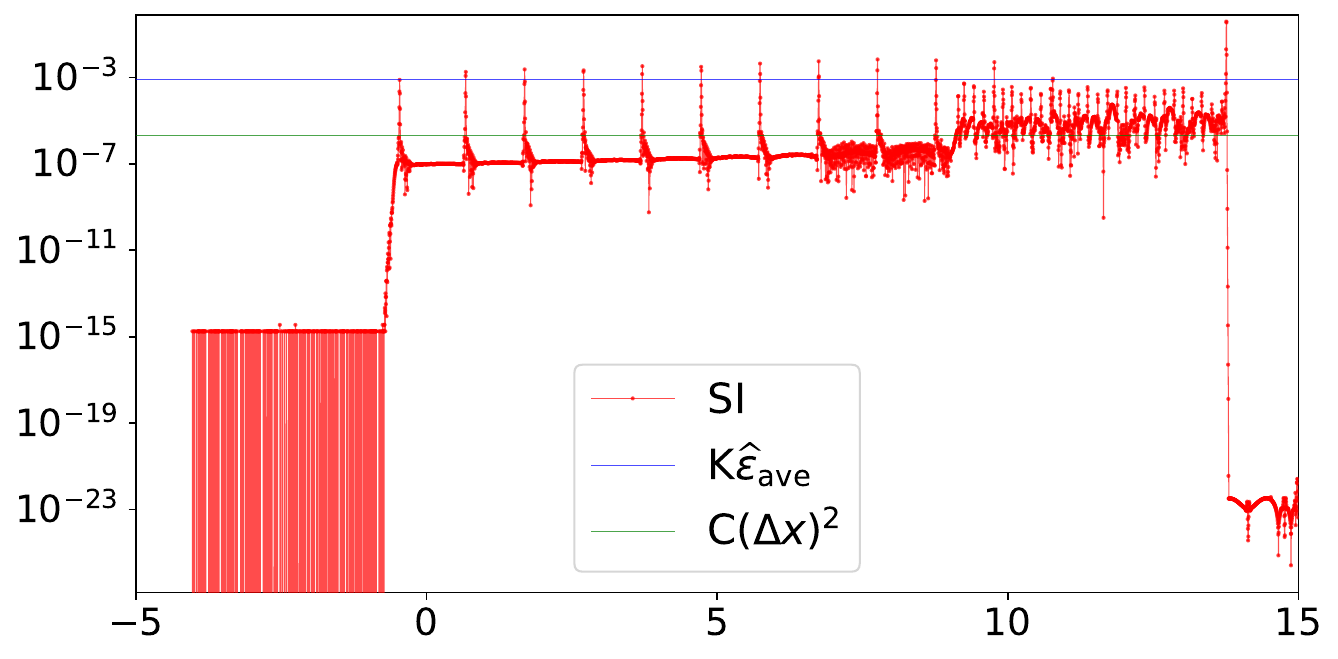}\hspace*{0.5cm}
            \includegraphics[trim=0.3cm 0.3cm 0.2cm 0.2cm, clip, width=5.4cm]{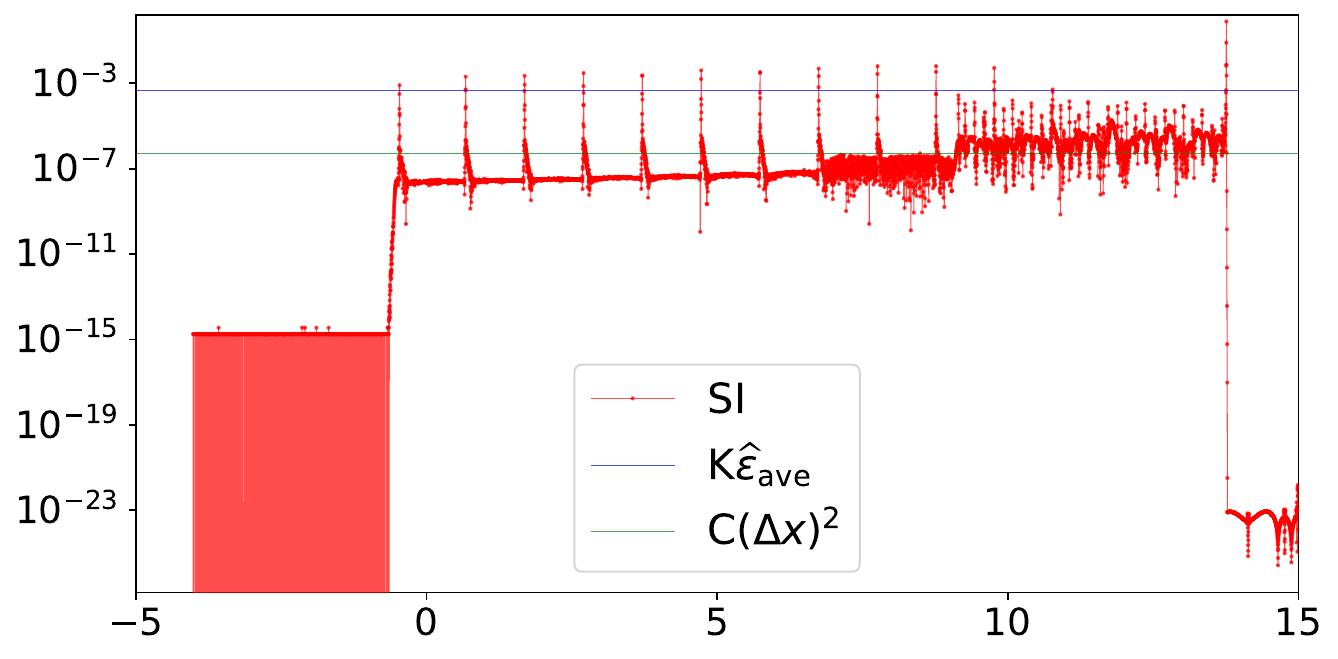}}
\caption{Example 2: Density computed using $N=800$ (top panel) and the corresponding SI values obtained on several meshes with $N=800$
(middle left), $N=1600$ (middle right), $N=6400$ (low left), and $N=12800$ (low right). $C=0.2$, $K=6$.\label{fig2}}
\end{figure}

\vskip5pt
\noindent
{\bf Example 2---shock-density wave interaction problem.}
The density computed with $N=800$ is plotted in Fig. \ref{fig2} (top panel). The corresponding SI values for various mesh refinements are
presented in the middle and low rows of Fig. \ref{fig2}. In this test, several shocks are separated by smooth areas and as one can clearly
see, the proposed SI successfully identifies them especially on a finer mesh, at which smaller shocks are more accurately captured by the
underlying AF method. As anticipated, the size of the SI decays proportionally to $(\dx)^2$ in the smooth areas.

\begin{figure}[ht!]
\centerline{\includegraphics[trim=0.2cm 0.3cm 0.2cm 0.3cm, clip, width=5.6cm]{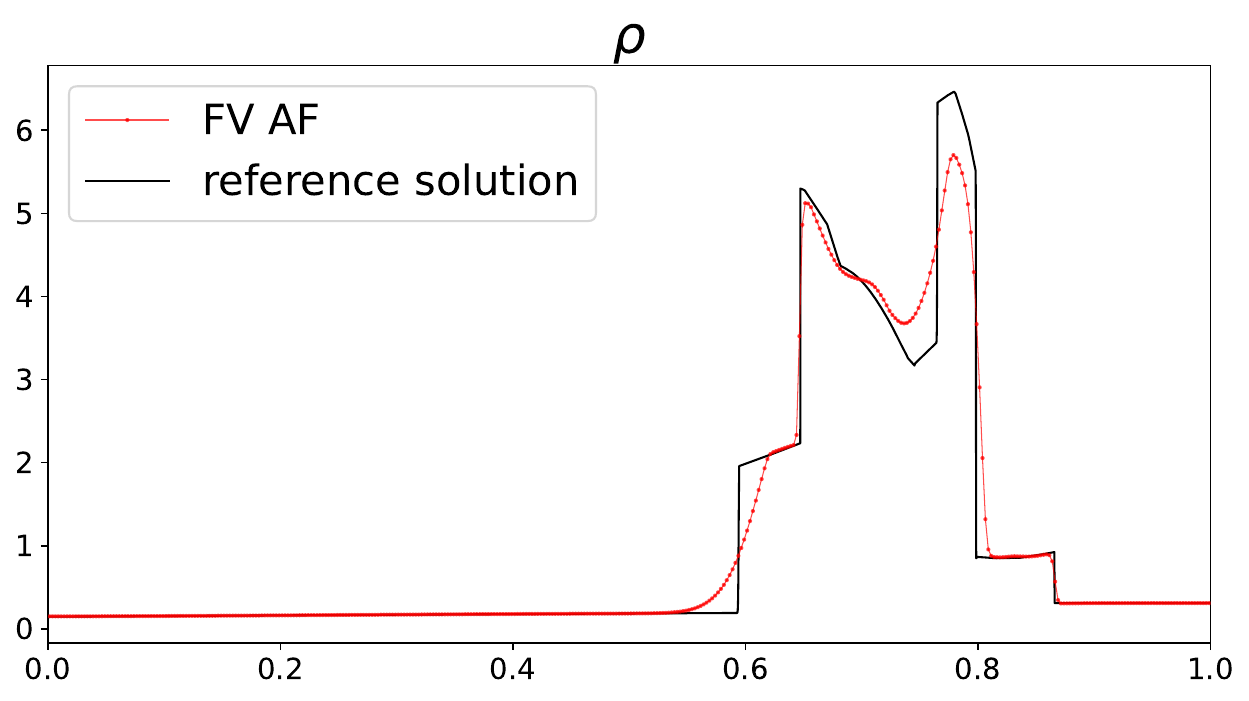}}
\vskip5pt
\centerline{\includegraphics[trim=0.3cm 0.3cm 0.2cm 0.2cm, clip, width=5.4cm]{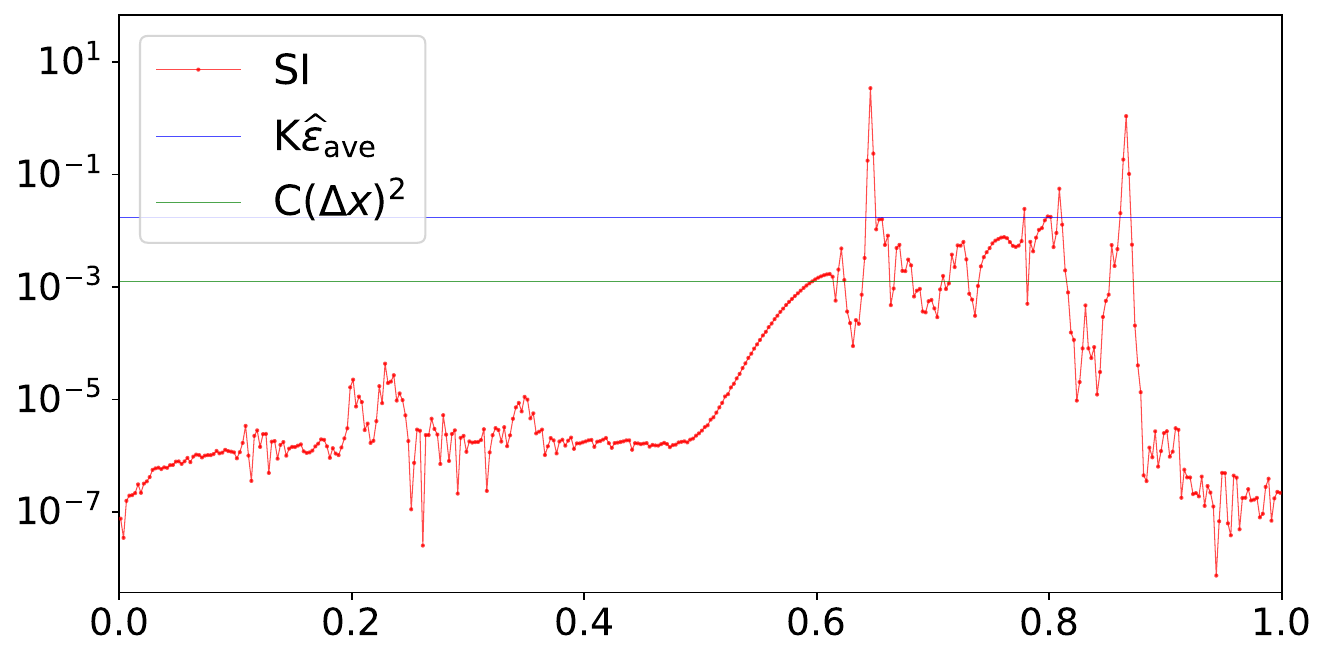}\hspace*{0.5cm}
            \includegraphics[trim=0.3cm 0.3cm 0.2cm 0.2cm, clip, width=5.4cm]{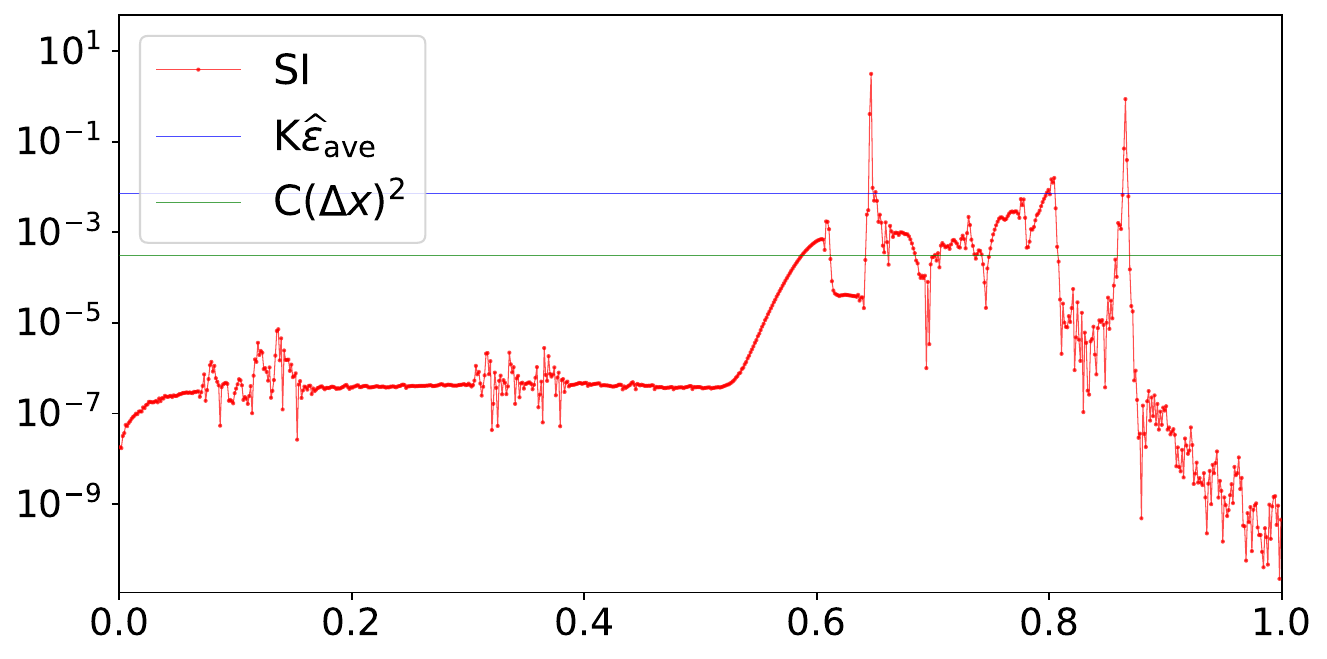}}
\vskip5pt
\centerline{\includegraphics[trim=0.3cm 0.3cm 0.2cm 0.2cm, clip, width=5.4cm]{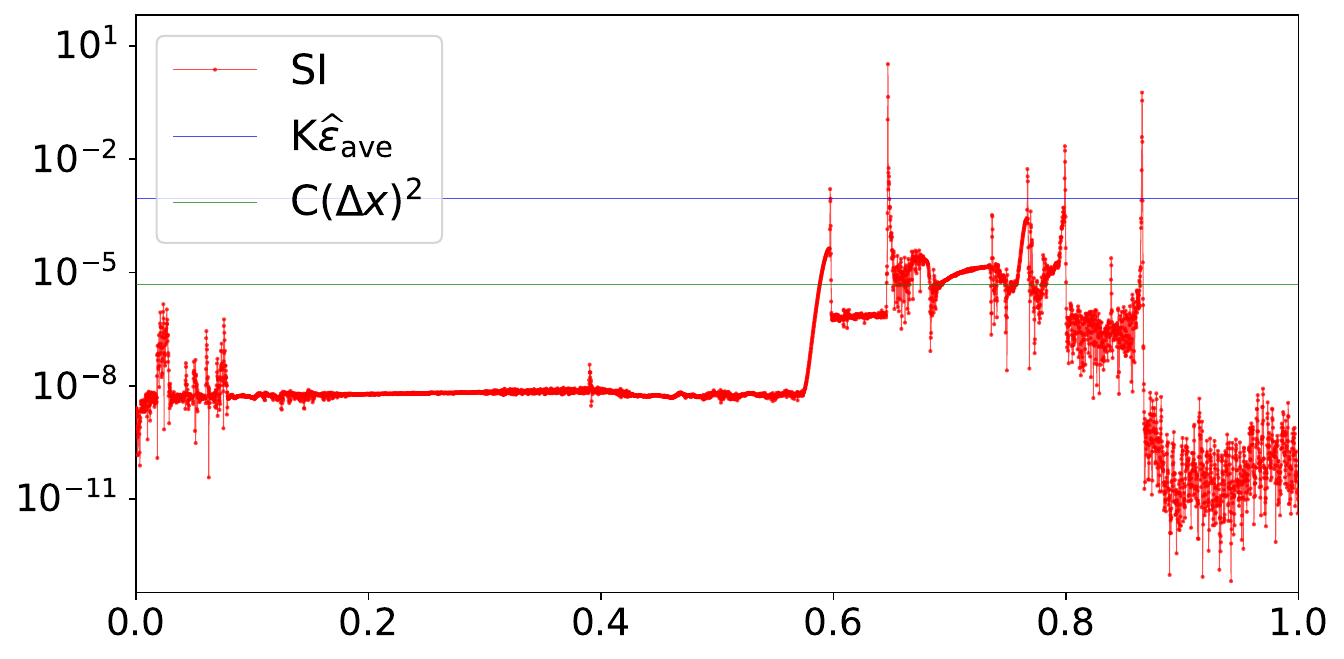}\hspace*{0.5cm}
            \includegraphics[trim=0.3cm 0.3cm 0.2cm 0.2cm, clip, width=5.4cm]{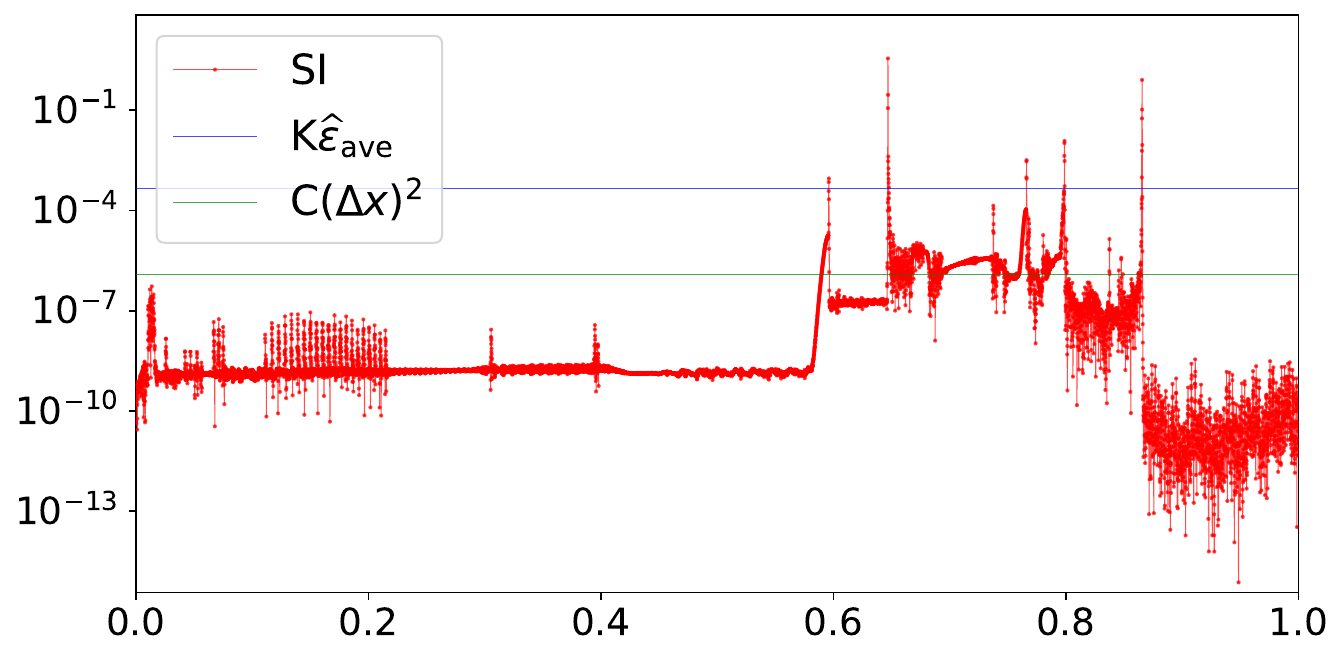}}
\caption{Example 3: Density computed using $N=400$ (top panel) and the corresponding SI values obtained on several meshes with $N=400$
(middle left), $N=800$ (middle right), $N=6400$ (low left), and $N=12800$ (low right). $C=200$, $K=1.2$.\label{fig3}}
\end{figure}

\vskip5pt
\noindent
{\bf Example 3---Woodward-Colella problem.}
Fig. \ref{fig3} shows the density computed with $N=400$. Additionally, the same figure presents SIs for various levels of mesh refinement.
Notably, even with different levels of refinement, the method consistently captures all key discontinuities in the solution, including the
left contact discontinuity, which is sharply captured on finer meshes only (see the bottom row in Fig. \ref{fig3}). This demonstrates the
robustness of the proposed SI.

\section{Concluding Remarks}
In this work, we have introduced a novel SI, which is based on measuring the difference between two solutions produced in the process of
solving hyperbolic (systems of) conservation laws by AF methods. While we have presented the new SI in the context of recently proposed
second-order semi-discrete FV AF method from \cite{ACKM}, it can be computed based on the solutions obtained by other AF methods as well as
other methods, in which two numerical solutions are evolved in time, for example, the methods on overlapping cells in
\cite{Liu2005,LSTZ2,LSX,XLDLS}.

The robustness of the proposed SI hinges on the expectation (confirmed by the presented numerical experiments) that in the ``rough'' parts
of the computed solutions the SI values are ${\cal O}(1)$, while in the smooth areas, its values decay as $(\dx)^2$ for the studied
second-order semi-discrete FV AF method. The developed SI is expected to be even more robust for higher-order methods as it will be able to
easier distinguish between ${\cal O}(1)$ and ${\cal O}((\dx)^r)$, where $r$ is the order of the method, and thus to better detect shocks and
especially contact discontinuities on coarser meshes.

We note that the multidimensional extension of the proposed SI are possible on both Cartesian and triangular (polygonal) meshes.

Finally, we emphasize that the new SI can be used to design a variety of adaptive strategies including a simple scheme adaption (in which,
for instance, different limiters are used in different parts of the computed solution) and adaptive mesh refinement.

\subsection*{Acknowledgment}
The work of A. Chertock was supported in part by NSF grant DMS-2208438. The work of A. Kurganov was supported in part by NSFC grants
12171226 and W2431004. The work of L. Micalizzi was supported in part by the LeRoy B. Martin, Jr. Distinguished Professorship Foundation.

\bibliographystyle{ws-procs9x6} 
\bibliography{sn-bibliography}

\end{document}